\documentclass[12pt]{amsart}

\usepackage{latexsym}
\usepackage{amsmath}
\usepackage{amssymb}
\usepackage{amscd}
\usepackage{psfig}
\usepackage{multicol}

\usepackage{latexsym, graphics, psfrag, amscd, amssymb, graphicx, psfrag}

\begin{document}

\author{Nan-Kuo Ho}
\address{Department of Mathematics\\ National Cheng-Kung University
 \\Taiwan}\thanks{The first author
was partially
supported by  grants from OGS and OGSST}
\email{nho@fields.utoronto.ca; nankuo@mail.ncku.edu.tw}

\author{Lisa  C. Jeffrey}
\address{Department of Mathematics\\ University of Toronto
 \\Toronto, Ontario M5S 3G3, Canada}
 \email{jeffrey@math.toronto.edu}\thanks{The second
author was partially supported by
a grant from NSERC}

\begin{abstract}
We compute the Riemannian volume of the moduli space of flat
connections on a nonorientable 2-manifold, for a natural
class of metrics. We also show that Witten's volume formula for
these moduli spaces may be derived using Haar measure, and we give a new
proof of Witten's volume formula for the moduli space of flat
connections on a 2-manifold using Haar measure.
\end{abstract}
\title[The Volume of the Moduli Space]{The
volume of the moduli space of  flat connections on a  nonorientable 2-manifold}

\maketitle

\newtheorem{thm}{Theorem}
\newtheorem{tm}[thm]{Theorem}
\newtheorem{ex}{Example}
\newtheorem{fact}[thm]{Fact}
\newtheorem{lm}[thm]{Lemma}
\newtheorem{rem}[thm]{Remark}
\newtheorem{cor}[thm]{Corollary}
\newtheorem{pro}[thm]{Proposition}
\newtheorem{de}[thm]{Definition}

\newcommand{\ra}{\rightarrow}
\newcommand{\Si}{\Sigma}
\newcommand{\tS}{\tilde{\Sigma}}
\newcommand{\M}{\mathcal{M}}
\newcommand{\Z}{\mathbb{Z}}
\newcommand{\R}{\mathbb{R}}
\newcommand{\hol}{\mathrm{Hom}(\pi_1(\Si),G)}
\newcommand{\ep}{\epsilon}

\newcommand{\A}{\mathcal A}
\newcommand{\cG}{\mathcal G}
\newcommand{\fG}{\mathfrak G}
\newcommand{\g}{\mathfrak g}
\newcommand{\tri}{\triangle}
\newcommand{\mf}{\mbox{flat}}
\newcommand{\mbas}{\mbox{bas}}
\newcommand{\Tra}{ {\rm Tr}}
\newcommand{\dimmalpha}{ { {\rm dim} (\alpha) }}

\newcommand{\scriptdee}{ {\mathcal{D}}}

\newcommand{\matr}[4]{\left \lbrack \begin{array}{cc} #1 & #2 \\
     #3 & #4 \end{array} \right \rbrack}

\newcommand{\gtt}{{g_{22}}}


\section{Introduction}

In \cite{w} Witten defined and computed a
volume on the moduli space of gauge equivalence classes
of flat connections on a
2-manifold, using
Reidemeister-Ray-Singer torsion (see e.g. \cite{Fr}).
 When the 2-manifold is orientable, Witten proved
that this volume is equal to the symplectic volume on the moduli space.
However, when the 2-manifold is not orientable the moduli space does not
have a symplectic form, so the interpretation of this volume has
been unclear.

The moduli space of gauge equivalence classes of flat connections on
a nonorientable 2-manifold turns out to be a Lagrangian submanifold
of the moduli space of  gauge equivalence classes of flat connections
on the orientable double cover.
In this article we compute its Riemannian volume for a natural class
of metrics on the 2-manifold. The dependence on the choice of 
metric is discussed in Remarks 1 and 2 below.

The layout of this article is as follows. In Section 2 we
compute the Riemannian volume of the moduli space of flat connections on
a nonorientable 2-manifold (exhibited as the connected sum of a
Riemann surface with either one or two copies of
$RP^2$) using a metric derived from a particular class of
metrics on $RP^2-  \{\rm disc\}$.
In Section 3
we give a new proof that Witten's formula for the volume of the moduli
space  of  flat connections on  2-manifolds (for nonorientable surfaces
with no boundary and orientable surfaces with one boundary component)
arises from Haar measure.

{\em Acknowledgement:} This article comprises
part of the Ph.D. thesis
of the first author, under the supervision of the second author.
Both authors would like to thank Eckhard Meinrenken  and
Chris Woodward for
useful conversations and insightful advice.

\section{Metrics on the moduli space}
\subsection{Preliminaries}
By the classification of 2-manifolds \cite{ma}, all nonorientable 2-manifolds
are obtained as the connected sum of a Riemann surface with either
one copy of 
the real  projective plane $RP^2$ (denoted $P$) or
two copies of $P$ (which is equivalent to the connected sum with 
a Klein bottle $K$).


Given a Riemannian metric on the M\"obius strip $P \setminus D$,
 we obtain a Riemannian metric on the connected sum 
$\Sigma \# P$ or $\Sigma \# P \# P$, where $\Sigma$ is 
an oriented 2-manifold equipped with a Riemannian metric. 
We assume that a Riemannian metric has
been chosen on $\Sigma \setminus D$.

We view $P \setminus D$ as formed by gluing together 
two sides of equal length of a  triangle $\tri$. 
We define a metric on $P \setminus D$ by endowing
$\tri$ with a Riemannian metric in which two sides are 
of equal length and  using geodesic 
polar coordinates about one vertex  of the  triangle.
In such coordinates $(\rho, \sigma)$ 
(where $\rho$ is the arc length from $p$
and $\sigma $ is the angular variable) the metric takes the form
(\cite{dC} Section 4.6)

\begin{equation} 
ds^2 = d\rho^2 + \gtt(\rho, \sigma) d\sigma^2 \end{equation} 

We now assume that the triangle $\tri$ is a geodesic triangle.
In geodesic polar coordinates, the Hodge star operator is then
\begin{equation} 
\ast d\rho = \sqrt{\gtt} d\sigma
\end{equation}
$$ \ast  d\sigma = - \frac{1}{\sqrt{\gtt}} d\rho $$
We make the assumption that 
$\gtt =\gtt(\rho)$ depends only on $\rho$ and is independent
of $\sigma$.
This assumption is satisfied by the three important 
examples of metrics with constant scalar curvature:
\begin{enumerate} \label{d:constcurv} 
\item Spherical metric
with scalar curvature $+1$:  $\gtt(\rho) = \sin^2 (\rho)$
\item Euclidean metric (with scalar curvature $0$):
$\gtt(\rho) = \rho^2$
\item Hyperbolic metric with scalar curvature $-1$:
$\gtt(\rho) = \sinh^2(\rho)$
\end{enumerate}

We can consider two triangles shown in Figure 1.
We use geodesic polar coordinates at the vertex $x_0$. 
The curves $\gamma_1$ and $\gamma_2$ are radial geodesics
with constant $\sigma$.
The curve $\gamma_3$ is the geodesic joining the two vertices
$x_1$ and $x_2$.
The curve $C_1$ is the curve $\rho $ = constant  between 
$x_1$ and $x_2$; it is not a geodesic.
The  triangle $\tri$ is bounded by 
$\gamma_1, \gamma_2, \gamma_3$; it is a geodesic triangle.
 The triangle 
$\tri_0$ is bounded by $\gamma_1, \gamma_2, C_1$.
We can find the integral over the triangle $\tri_0$
bounded by $\gamma_1,\gamma_2,C_1$
where $C_1$ is the
curve specified by the equation  $\rho=L =$ constant.
\begin{figure}[!htbp]
        \begin{center}
                \psfrag{r1}[c][c][1][0]{$\gamma_1$}
                \psfrag{r2}[c][c][1][0]{$\gamma_2$}
                \psfrag{r3}[c][c][1][0]{$\gamma_3$}
                \psfrag{d}[c][c][1][0]{${\mathcal D}$}
                \psfrag{p}[c][c][1][0]{$x_0 \phantom{aa}\phi $}
                \psfrag{c1}[c][c][1][0]{$C_1$}
                \includegraphics[scale = 0.7]
                {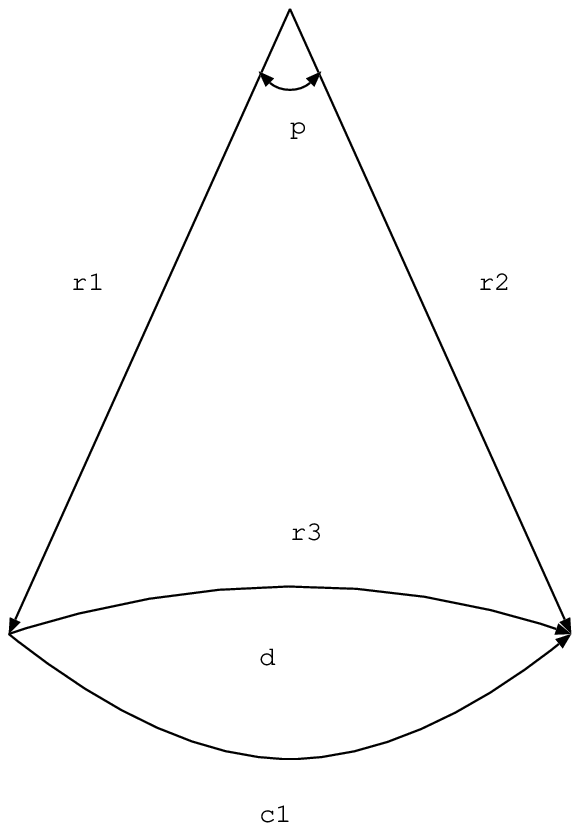}
        \end{center}
       \caption{Triangles $\tri$ and $\tri_0$}
\end{figure}

Let $G$ be a compact connected Lie group;
the moduli space of flat connections on
$\tri$ modulo based gauge transformations (i.e. gauge transformations
which take the value identity on the three vertices) can be identified with
$G \times G$.
The moduli space of flat connections on $P\setminus\scriptdee$ modulo
gauge transformations which are trivial at one point on the boundary
(corresponding to the vertices of the geodesic triangle) is
equivalent to the subspace
of the  moduli space of flat connections on the  geodesic
triangle (modulo based gauge transformations) which correspond
under an orientation reversing  map    which identifies two sides of the triangle .
This is ${\rm Hom}(\pi,G)$ where $\pi = \pi_1(P \setminus D)$,
in other words $${\rm Hom}(\pi, G)
 = \{ (u,v) \in G \times G ~|~ u v^2= 1 \} \cong G.$$
Here $v$ is the holonomy along $\gamma_1$ and $- \gamma_2$,
while
$u$ is the holonomy along $\gamma_3$.

On the space of flat connections on $\tri$, there is a metric
\begin{equation}
<a,b>=\int_\tri\Tra(a\wedge\ast b)\end{equation}
where $\ast$ is the Hodge star operator on differential forms over $\tri$.
We will study  the corresponding metric on $G \times
G$. Define the metric $<,>$ at a flat connection $A$
corresponding to a point $p$ in $G \times G$ to be
\[<d_A\xi,d_A\eta>=:\int_{\tri}\Tra(d_{A}\xi\wedge\ast
d_{A}\eta)\] where $\xi(\cdot)$ and $\eta(\cdot)$ are $\g$-valued 
functions on $\tri$, so $d_A \xi$ and $d_A \eta$ represent
elements of the tangent space to $G\times G$ at $p$.
Here,  $\Tra(\cdot)$ represents the ad-invariant inner
product on the Lie algebra.

In the case when $G$ is abelian, the holonomy of $d \xi$ along $\gamma_1$ is
\begin{equation} \label{holone} Hol_{\gamma_1}d\xi=\exp\int_{\gamma_1}d\xi=\exp(\xi(\rho=L_1,\sigma=0)-\xi(\rho=0,\sigma=0)), \end{equation}
and its holonomy along $\gamma_2$ is
\begin{equation} \label{holtwo} Hol_{\gamma_2}d\xi=\exp\int_{\gamma_2}d\xi=\exp(\xi(\rho=L_2,\sigma=0)-\xi(\rho=0,\sigma=\phi)). \end{equation}
We will study these holonomies once we have determined $\xi$.

From now on, we only
consider the inner product at $A=0$. The reason is that the metric
$<a,b>=\int \Tra(a\wedge \ast b)$ is ad-invariant
(where $\ast$ denotes the Hodge star operator),
and since the
gauge action at the infinitesimal level is the adjoint action, this
metric is invariant under the action of the gauge group and in
particular under the action of the based
gauge group. Since 
the space $\tri$ is contractible, 
any flat connection is gauge equivalent to the
trivial connection on a triangle.
In other words all infinitesimal flat connections -- elements
 of the
tangent space of the space of flat connections modulo based gauge transformations -- can
be expressed as $d\xi$ for some infinitesimal gauge transformation $\xi:\tri\ra\g$. So we need only  consider the
inner product at $A=0$.

We want to associate a norm 
$< d\xi, d\xi>   $ where
$\xi \in \Omega^0(P \setminus D) \otimes \g$. 
This will be done by defining

\begin{equation} \label{e:norm}
<d\xi,d\xi>=:\int_{\tri}\Tra(d\xi\wedge\ast
d\xi)\end{equation}
Using Stokes' theorem,
 this can be converted to a line integral around the boundary of
the triangle if we have $d\ast d\xi=0$: this follows because
\begin{eqnarray*}
\int_{\tri}\Tra(d\xi\wedge\ast
d\xi)&=&\int_{\tri}\Tra[d(\xi\wedge\ast
d\xi)-\xi d\ast d\xi]\\
&=&\int_{\partial\tri}\Tra\xi \ast d\xi -\int_{\tri}\Tra\xi
d\ast d\xi
\end{eqnarray*}

In fact, we may assume $\xi$ satisfies
\begin{equation}\label{e:3.2} d\ast d\xi=0.\end{equation}
This condition represents a transversal to the orbit of the based gauge group; our
procedure is analogous to identifying
 de Rham cohomology classes with harmonic forms
(forms $\alpha$ satisfying $d\alpha=0$ and $d\ast\alpha=0$). Note that the
space $G\times G$ is isomorphic to the space of 
flat connections on $\tri$ modulo based
gauge transformations (i.e. gauge transformations 
which take the value $1$ at the vertices).
Each equivalence class may be written as $d\xi$ for 
some $\xi:\tri\ra \g$ (which does
not take value $1$ at all the vertices, 
unless one wishes to represent the trivial
flat connection).

To solve equation (\ref{e:3.2}), first, let us recall some
 geometry.
Recall that our 
Riemannian metric in geodesic polar coordinates is
\[ds^2=d\rho^2+\gtt(\rho)  d\sigma^2\]where $\rho$ is the
distance from a chosen point $x_0$ 
 and $\sigma$ is the polar angle
with respect to   this point.
Let $\ast$ denote the Hodge star operator on the 2-manifold;
the star operator with respect to
 these polar coordinates
is\[ \ast d\rho=\sqrt{\gtt(\rho)}
 d\sigma, ~~\ast d\sigma=\frac{-1}{\sqrt{\gtt(\rho)}} d\rho. \]

 We use the
substitution 
\begin{equation} \label{e:udef}
u(\rho)=\int^{\rho}\frac{1} { \sqrt{\gtt( t)} } dt 
\end{equation}
and define
$$\tau(\rho) = \exp u(\rho).$$
The equations for the Hodge star operator become
\begin{equation} \label{hodge2}
* du = d \sigma \end{equation}
$$ * d\sigma = - du $$
We now have to
solve the equation
$$ d (* d \xi) = 0 $$
We have
$$ d\xi = \frac{\partial \xi}{\partial u} du +\frac{\partial  \xi}
{\partial\sigma} d\sigma $$
so
\begin{eqnarray*} d (* d \xi) &=&d( \frac{\partial  \xi}{\partial u} d\sigma - \frac
{\partial \xi}{\partial \sigma}  du)
  =( \frac{\partial^2  \xi}{\partial \sigma^2} +
\frac{\partial^2  \xi}{\partial u^2} ) du \wedge d\sigma
\\&=&(\xi_{\sigma\sigma}+\xi_{uu})du \wedge d\sigma
\end{eqnarray*}

So the equation we need to solve is
\begin{equation} \label{eq:secorder}
 \xi_{\sigma\sigma}+\xi_{uu}  = 0. \end{equation}

Let us take a Fourier decomposition of $\xi$ in the $\sigma$ variable:
$$\xi(u,\sigma)=\int_{\R} \hat{\xi}(u,k)e^{ik\sigma}dk $$
So (\ref{eq:secorder}) becomes
\begin{equation} \label{eq:consequence}
\frac{\partial^2}{\partial u^2} \hat{\xi} (u,k) - k^2 \hat{\xi} (u,k) = 0
\end{equation}

This equation (\ref{eq:consequence}) has the following solutions:

$$\hat{\xi} (u,k) = C_+ (k) \tau^k + C_-(k) \tau^{-k}. $$
Recall that  we had defined  $\tau = \exp (u)$.
Imposing the condition that $\hat{\xi} (u,k)$ is finite at
$u = 0 $, we get
$$ \hat{\xi} (u,k) = C_+ (k) \tau^k $$
when $k > 0 $,  and
$$ \hat{\xi} (u,k) = C_- (k) \tau^{|k|} $$
when $k < 0$.
Recall that $\phi$ is the polar angle of the geodesic triangle.
Let us impose the additional constraint that
$\xi(u, \sigma)$ attains its maximum on two edges of the
 triangle ($\gamma_1$ which is
$\sigma = 0 $ and $\gamma_2$ which is  $\sigma = \phi$),
in other words
$\frac{\partial}{\partial \sigma} \xi(\rho, \sigma) = 0 $
when $\sigma = 0 $ and $\sigma = \phi$. This leads to
$$ \xi(u, \sigma) = Y\tau^k \cos k \sigma $$
where $Y \in \g$ is a constant and
\begin{equation} \label{e:kdef} k = \pi/\phi. \end{equation}

We compute that 
\begin{equation}\int_{C_1} \Tra (\xi \wedge \ast d\xi)  
= \int_0^\phi <\xi, \partial \xi/\partial u > d\sigma \end{equation}
$$ = < Y, Y> k \exp k u(L) \int_0^\phi \cos^2 k \sigma d\sigma $$
$$ = <Y,Y> \pi/2 \exp k u(L) . $$

Note that writing $\xi$ in terms of its Fourier decomposition
accomplishes the same thing as solving the equation $d * d \xi = 0 $
by the method of separation of variables.
We solved this equation in the previous subsection.
Thus we obtain \begin{equation} \label{e:xi}
\xi(\rho,\sigma)=Y\tau(\rho)^k\cos k\sigma,~~~~k=\pi/\phi
\end{equation}
so that the maximum of $|\xi|$ is achieved on the edges $\sigma = 0 $
and $\sigma = \phi$
and $Y \in\g$ is a variable we choose so that we can get the desired
holonomy along the boundary.
The holonomies were specified by the equations
(\ref{holone}) and (\ref{holtwo}).

\begin{equation} \label{e:gamone}
\int_{C_1}\Tra(\xi\wedge\ast d\xi)= 
<Y,Y>k\Bigl [\tau (\frac{L}{2})\Bigr ]^{2k}\int_0^\phi d\sigma
\cos^2k\sigma \end{equation}
$$
=<Y,Y>\Bigl [ \tau (\frac{L}{2})\Bigr ]^{2k}\frac{\pi}{2} $$
$$ \int_{\gamma_1}\Tra(\xi\wedge\ast d\xi)= 
\int_{\gamma_2}\Tra(\xi\wedge\ast d\xi)=0
$$

Thus (referring to Figure 1)
the integral over the geodesic triangle $\tri$ bounded by
$\gamma_1,\gamma_2,\gamma_3$ is the integral over the triangle
$\tri_0$ 
bounded by $\gamma_1,\gamma_2,C_1$ minus the integral over the
region ${\scriptdee}$:
\begin{equation} \label{e:components}
\int_\tri \Tra(d\xi\wedge\ast d\xi) = 
\int_{\tri_0} \Tra(d\xi\wedge\ast d\xi)  - 
\int_{\scriptdee} \Tra(d\xi\wedge\ast d\xi).
\end{equation}
This leads to 
\begin{equation} \label{e:tri}
\int_\tri  \Tra(d\xi\wedge\ast d\xi) =  \int_{C_1}
  \Tra(\xi\wedge\ast d\xi) - \int_{\scriptdee}  \Tra(d\xi\wedge\ast d\xi) .
\end{equation}
The integral over $\scriptdee$ is given as follows:

$$
\int_{{\scriptdee}}\Tra(d\xi\wedge\ast d\xi) =
\int_{{\scriptdee}}<Y,Y>k^2\tau^{2k}
(\frac{\cos^2k\sigma}{\tau}+\frac{\sin^2k\sigma}{\tau}
)d\tau\wedge d\sigma$$
$$ = \int_{{\scriptdee}}<Y,Y>
k^2\tau^{2k}\frac{1}{\tau}d\tau\wedge d\sigma $$
$$ =<Y,Y>k^2\int_0^\phi d\sigma
\int_{\tau(\gamma_3(\sigma))}^{\tau(C_1(\sigma))}\tau^{2k}
\frac{d\tau}{\tau} $$
\begin{equation} \label{tangamdef}
= <Y,Y>k^2\int_0^\phi (\frac{\tau^{2k}(C_1(\sigma))}{2k}-
\frac{\tau^{2k}(\gamma_3(\sigma))}{2k})d\sigma
\end{equation}
where $\tau(C_1(\sigma))$ is a constant independent of $\sigma$.

\begin{rem} If $g$ is a Riemannian metric on $\tri$ then the quantity
$$<d\xi, d\xi>_g $$ (associated to the metric $g$)
 is equal to the 
quantity 
 $<d\xi, d\xi>_{ge^w} $ 
where $w: P \setminus D \to \R$ is a 
$C^\infty$ function. In other words our definition 
of the metric  $<d\xi,d\xi>_g$ is invariant under
conformal transformations on $\tri$. 
To see this, we 
write the metric as 
$$ g = \matr{g_{11}}{g_{12}}{g_{21}}{g_{22}}$$
with inverse
$$ g^{-1} = \matr{g^{11}}{g^{12}}{g^{21}}{g^{22}}.$$
We
observe that the Hodge star operator
is unchanged under $g \mapsto g e^w$ (where $w: \tri \to \R$),
since
$\det g$ transforms to $(\det g) e^{2w}$ while
$g^{11}$ transforms to $g^{11} e^{-w}.$ 
The Hodge star operator for  coordinates $x_1, x_2$ is 
$$\ast dx_1 = {g^{11}}  \sqrt{\det g}dx_2 $$
Hence $g^{11}\sqrt{\det g} $ is unchanged. 
Thus the norm (denoted by $\Gamma$)  on $d\xi$ computed using
metrics of the form 
\begin{equation}
\label{e:form}  ds^2 = d\rho^2 + \gtt(\rho) d\sigma^2 \end{equation}
also gives the norm for metrics $g$
conformally equivalent to those of the form (\ref{e:form}).
\end{rem}
\begin{rem} In fact in dimension 2 every Riemannian metric is 
locally diffeomorphic to one which is conformally 
equivalent to a metric of constant curvature
(see \cite{d'H} Section 3.3). 
This shows that 
all metrics on $P \setminus D$ are conformally equivalent to 
a metric of constant curvature, one of the three listed in  
(\ref{d:constcurv}), for which the 
norm (\ref{e:norm}) may be computed as in (\ref{e:components}).
The values of the norms are  different for the three
different choices of constant curvature metrics. 
\end{rem}

\begin{rem}
The volume of the moduli space defined using Reidemeister-Ray-Singer torsion
is independent of the choice of metric on the 2-manifold either orientable or
nonorientable. The Riemannian volume, however, is a different story. For an
orientable 2-manifold, the moduli space of flat connections is a K\"ahler manifold,
so the Riemannian volume equals the symplectic volume, and is thus independent of
the choice of metric. For a nonorientable surface, Witten remarks
in \cite{w} (p. 163) that the Riemannian volume does depends on the metric.
According to \cite{w} (2.38)
the torsion volume and the Riemannian volume differ (in this case) by a ratio
of determinants of elliptic operators
(generalized Laplacians on Lie algebra-valued differential
forms).
Our results are consistent with this observation, and at the same time provide
an evaluation of this ratio of determinants.
\end{rem}

\subsection{Hyperbolic metric on the geodesic triangle}
We explore the hyperbolic
case first. Note that if $\Sigma$ has genus $\ge 2$ then by 
the uniformization theorem it is equipped with a unique metric
of constant scalar curvature $-1$ (in other words a hyperbolic
metric). This is one motivation for choosing
a hyperbolic metric on $P \setminus D$, although 
such a metric will be singular at one point. 




A Riemann surface with one boundary component and genus $\ge 1$
always has a hyperbolic metric with constant curvature
$ -1$ \cite{bu}. For this reason we construct a
metric on $P \setminus D$ 
for which the
boundary is a geodesic and which is hyperbolic with constant curvature
$-1$ at all but one point. We do this by gluing the edges of an isosceles
geodesic hyperbolic triangle using an orientation reversing 
map:
this process identifies all three vertices
of the triangle, and the vertex becomes the point where the metric does
not have constant curvature $-1$ (indeed, the metric is singular at this
point).

\begin{rem}\label{singular}
The reason why this procedure gives a natural metric on
the 2-manifold formed by taking the connected sum
of a Riemann surface of
 genus $\ell>0$ with
$P$ is that a Riemann surface of genus $\ell > 0 $
with one or two boundary components always
has a hyperbolic metric with constant curvature for which the boundary
components are geodesics. This metric can be obtained from a pants
decomposition of the Riemann surface, where each pair of pants is equipped
with a hyperbolic metric for which the boundary 
components are geodesics. See \cite{bu}.
\end{rem}

Hence the connected sum of a Riemann surface with $P$ has a metric which
is hyperbolic at all but one point, and the connected sum of a Riemann
surface with two copies of $P$ has a metric which is hyperbolic at
all but 2 points.
To form  $P$ with one disc removed, we recall that the two edges
of the triangle  are identified using an
orientation reversing map. In Figure 1 this corresponds to identifying
$\gamma_1$ with $\gamma_2$ using an orientation reversing map which maps
$x_1 \in \gamma_1$ to $x_0 \in \gamma_2$ and maps $x_0 \in \gamma_1$ to
$x_2 \in \gamma_2$. It turns out that the metric on the based moduli
space on $P$ with one disc removed is singular at the vertex.

\begin{figure}[!htbp]
        \begin{center}
                \psfrag{x0}[c][c][1][0]{$x_0$}
                \psfrag{x1}[c][c][1][0]{$x_1$}
                \psfrag{x2}[c][c][1][0]{$x_2$}
                \psfrag{r1}[c][c][1][0]{$\gamma_1$}
                \psfrag{r2}[c][c][1][0]{$\gamma_2$}
                \psfrag{r3}[c][c][1][0]{$\gamma_3$}
                \psfrag{ph}[c][c][1][0]{$\phi$}
                \includegraphics[scale = 1.5]    {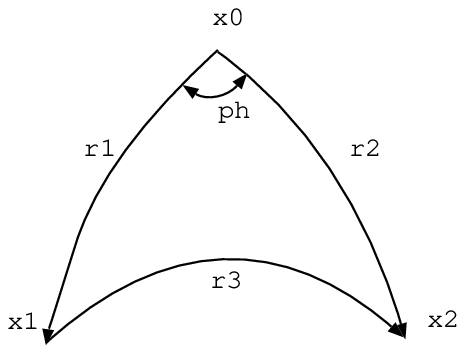}
        \end{center}
        \caption{Hyperbolic geodesic triangle}
\end{figure}
We have introduced a geodesic triangle $\tri$ in the upper half plane
 (see Figure 2).
As shown in Figure 2, the geodesic triangle has one distinguished
vertex $x_0$  (which  serves as the origin for geodesic
polar coordinates). The sides $\gamma_1, \gamma_2, \gamma_3$ are
geodesics. We introduce geodesic polar coordinates $(\rho,\sigma)$
(using the hyperbolic metric on the upper half
plane, which is assumed to contain the
triangle) with $x_0$ as the origin:
 here $\rho$ is
 the distance from a chosen point $x_0$ and $\sigma$ is the polar angle.
It is assumed that the geodesic $\gamma_1$
is defined by the equation $\sigma = 0$,  and $\gamma_2$
is defined by $\sigma = \phi$ (where $\phi$ is the angle at
$x_0$). The curves $\gamma_1$ and
$\gamma_2$ have lengths $L_1$ and $L_2$ which are the same here.
The geodesic $\gamma_3 $  is specified
by Coxeter's equation
\begin{equation}\label{cox376}
\coth\rho=\frac{\cos(\sigma-\sigma_0)}{\ell}, ~~\mbox{ see
\cite{cox} p.376},\end{equation} where $\ell,\sigma_0$ are
constants determined by the geodesic $\gamma_3$. Also,
\[\coth\rho=\frac{\cosh\rho}{\sinh\rho}=\frac{\cosh^2(\rho/2)+\sinh^2(\rho/2)}{2\sinh(\rho/2)\cosh(\rho/2)}
=\frac{1}{2}(\tau+\frac{1}{\tau})\] for $\tau=\tanh(\rho/2)$.
Thus
\[
(1-\frac{\cos(\sigma-\sigma_0)}{\ell\tau})d\tau=\frac{-\sin(\sigma-\sigma_0)}{\ell}d\sigma
.\]
By (\ref{cox376}), it is clear that since $\rho (\sigma = 0 ) $
must be equal to $\rho(\sigma = \phi)$ for
a hyperbolic geodesic triangle two sides of which will be identified
to form $P \setminus D$, we need $\sigma_0 = \phi/2$ in (\ref{cox376}).

Now the equation for $\gamma_3$ is (see  equation
(\ref{cox376}))
\[ \cos(\sigma-\phi/2)\tanh\rho=\cos(\phi/2)\tanh L\] and
$$\tanh\rho=\frac{2\tau}{1+\tau^2}.$$
 So we have
\begin{equation} \label{taugdef}
\tau(\gamma_3(\sigma))=\frac{\cos(\sigma-\phi/2)}{l}
-\sqrt{(\frac{\cos(\sigma-\phi/2)}{l})^2-1},~~l=\cos(\phi/2)\tanh L
\end{equation}

Thus by (\ref{tangamdef})
\begin{eqnarray*}
 \int_{{\scriptdee}}\Tra(d\xi\wedge\ast d\xi)&=&
 \frac{k\phi}{2} <Y,Y>\Bigl [\tanh(\frac{L}{2})\Bigr ]^{2k}
\end{eqnarray*}
$$
-\frac{k}{2} <Y,Y>\int_0^\phi
\Bigl [\frac{\cos(\sigma-\phi/2)}{l}
-\sqrt{(\frac{\cos(\sigma-\phi/2)}{l})^2-1}~\Bigr ]^{2k}
d\sigma.$$
 Let
$x=\cos(\sigma-\phi/2)$;
\begin{eqnarray} \label{e:dint} \int_{{\scriptdee}}\Tra(d\xi\wedge\ast
d\xi) &=& \frac{k\phi}{2} <Y,Y>\Bigl [ \tanh(L/2)\Bigr ]^{2k}
\end{eqnarray}
$$   -\frac{k}{2} <Y,Y> \Bigl [
\int_{\cos(\phi/2)}^{1}\Bigl (\frac{x}{l}-\sqrt{(\frac{x}{l})^2-1}~
\Bigr )^{2k}
\frac{dx}{\sqrt{1-x^2}} $$
$$ -\int_1^{\cos(\phi/2)}
\Bigl (\frac{x}{l}-\sqrt{(\frac{x}{l})^2-1}~
\Bigr )^{2k}\frac{dx}{\sqrt{1-x^2}}\Bigr ]$$
$$ = <Y,Y>\frac{k\phi}{2} \Bigl [ \tanh (L/2)\Bigr ]^{2k}
-k <Y,Y> \int_{\cos(\phi/2)}^1\Bigl (\frac{x}{l}-\sqrt{(\frac{x}{l})^2-1}~
\Bigr )^{2k}
\frac{dx}{\sqrt{1-x^2}}
$$

Consider a more sp
ecific geodesic triangle, the triangle with two
sides of equal length $L$. This means $\sigma$ is fixed on those
two sides and the third side is determined by Coxeter's equation
(\ref{cox376}), where $(\rho,\sigma)$ are the polar coordinates as
defined above.
By (\ref{e:xi}) we have 
$$\xi(\rho, \sigma) =Y  (\exp u(\rho))^{k} \cos k \sigma $$
where 
in the notation of (\ref{e:udef}), $u = \tanh (\rho/2)$ in the case of
a hyperbolic metric.

We recall  from (\ref{e:tri}) that
$$
\int_{\tri}\Tra(d\xi\wedge\ast d\xi) =
\int_{C_1}\Tra(\xi\wedge\ast d\xi)-
\int_{{\scriptdee}}\Tra(d\xi\wedge\ast d\xi) $$
\begin{equation} \label{e:8a}
= k <Y,Y> \int_{\cos(\phi/2)}^1(\frac{x}{l}-\sqrt{(\frac{x}{l})^2-1}~)^{2k}
\frac{dx}{\sqrt{1-x^2}}
\end{equation}
We will compute this integral in Appendix A.

We denote the  integral (\ref{e:8a}) by $k <Y,Y> h(\phi)$ because
although it depends on two parameters $\phi$ and $L$, there
is a relation between these two.
 We can see this as follows. 
Assume the length of $\gamma_3$  is fixed (since $\gamma_3$ will
be glued to the boundary of a Riemann surface); denote this length by
$b$. We have

$$ 
\cosh^2L=\frac{\sinh^2b}{2(\cosh b-1)}+
\frac{(\cosh b-1)}{2}\frac{(1+\cos\phi)^2}{1-\cos^2\phi}
$$ 
$$ \frac{2(\cosh b-\cos\phi)(\cosh b-1)(\cos\phi+1)}
{2(\cosh b-1)(1-\cos^2\phi)}
$$ 
\begin{equation}\label{e:lphidep}
= \frac{\cosh b-\cos\phi}{1-\cos \phi}
\end{equation}
This is the relation between $\phi$ and $L$ once $b$ is fixed.

If $\tau(L)=\tanh L/2$, then the relation between $\tau(L)$ and $\phi$ is
\begin{equation}\label{tau}
\tau^2(\phi)=\tau^2(L(\phi))=
\frac{\sqrt{\frac{\cosh b-\cos\phi}{1-\cos\phi}}-1}{\sqrt{\frac{\cosh b-\cos\phi}{1-\cos\phi}}+1};
\end{equation}
thus $h$ is a function of the top angle $\phi$ only.

Thus from (\ref{e:8a}) we have
 \begin{equation}\label{rpmetric}
 <d\xi,d\xi>=k <Y,Y> h(\phi)
 \end{equation}

If we choose an orthonormal basis ${e_1,e_2,\cdots,e_n}$ for
$Lie(G)$, then the volume element of this moduli space using the
metric   $\Gamma$ is $\sqrt{\mbox{det}\Gamma}$, and
\[\mbox{det}\Gamma=(k h(\phi))^n,
~~~~\mbox{using}~~(\ref{rpmetric})\]
where $n$ is the dimension of $Lie(G)$ i.e. $\g$.

\subsection{Euclidean metric on the geodesic triangle}

The only smooth 
constant curvature metric on $P \setminus D$ for which the 
boundary is a geodesic  is the 
Euclidean metric.
For the geodesic triangle $\tri$ equipped 
with the Euclidean metric, we get 
$$u = \int^\rho \frac{d\rho}{\rho} = \ln \rho$$
so $\tau = \rho$ and $$\xi  = Y \rho^k \cos k \sigma.$$
The geodesic $\gamma_3$ is given by 
$$ \rho \cos (\sigma - \phi/2) = L $$
so we have
\begin{equation} 
\int_{C_1} \Tra (\xi \wedge \ast d \xi) = 
\frac{\pi}{2}<Y,Y> L^{2k} 
\end{equation}
and from (\ref{e:components})  and (\ref{tangamdef})
$$\int_{\scriptdee} \Tra (d\xi \wedge \ast d \xi) =   
\frac{\pi}{2} <Y,Y> L^{2k} - 
\frac{k}{2}<Y,Y> \int_{\rho \cos (\sigma - \phi/2) = L }  \rho^{2k} d\sigma 
$$
\begin{equation}=  
\frac{\pi}{2}<Y,Y> L^{2k} -k <Y,Y>  L^{2k} \int_0^{\phi/2}  
\sec^{2k}(\sigma) d\sigma . \end{equation}
Thus the norm $\int_{\scriptdee} \Tra (d\xi \wedge \ast d \xi) $
is given by 
$k <Y,Y> h (\phi)$ for an appropriate
function $h (\phi)$ which depends on the angle $\phi$ and 
on the choice of Euclidean metric.

\subsection{Spherical metric on the geodesic triangle}

For a spherical metric on the geodesic triangle,
we get 
\begin{equation} 
u = \int^\rho \frac{dx}{\sin x} = \ln \tan (\rho/2) \end{equation}
so $\tau = \tan (\rho/2).$
We then get that 
\begin{equation}
\int_{C_1}  \Tra (\xi \wedge \ast d \xi) = 
\frac{\pi}{2}<Y,Y> (\tan (L/2))^{2k} \end{equation}
Hence we get that
\begin{equation}
\int_{\scriptdee} \Tra (d\xi \wedge \ast d \xi) = 
<Y,Y>  \Biggl \{
\frac{\pi}{2}(\tan L/2)^{2k} - 
\frac{k}{2} \int_{\gamma_3} (\tan \rho/2)^{2k} d\sigma \Biggr \}. 
\end{equation}
Here, $\gamma_3$ is the segment of the great circle connecting the points
$x_1 = (\rho = L, \sigma = 0 )$ and 
$x_2 = (\rho = L, \sigma = \phi )$.
The length of this geodesic and the angles at
the vertices $x_1, x_2$ can be determined using spherical trigonometry
\cite{Weis}.
Thus 
when $\tri$ is equipped 
with a spherical metric of constant curvature,
the norm $<d\xi , d \xi> $
is given by 
$k <Y,Y> h (\phi)$ for an appropriate
function $h (\phi)$ which depends on $\phi$ and on the choice of spherical metric.

\subsection{$\Si$ is the connected sum of a Riemann surface with
$P$} \label{ss:csp}

To consider the volume of the moduli space of the connected sum of a Riemann
surface with $P$, first, let us examine the following lemmas:

\begin{lm}
Let $E$ be a Riemannian manifold with metric $\Gamma$.
Let $f:E\ra \R$ be a smooth
function for which $df_m \neq 0$ for any $m \in E$. Then
$$Vol(E)=\int_{\R}\frac{Vol f^{-1}(t)}{df(v(t))}dt$$
where $v(t)$ is the unit normal vector to $f^{-1}(t)$.
\end{lm}
\paragraph{\bf Proof}
The volume form on $T_mE$ is $e_1^{\ast}\wedge\cdots\wedge
e_n^{\ast}$ where $\{e_j\}$ is
an orthonormal basis of  tangent vectors i.e.
$\Gamma(e_j,e_k)=\delta_{jk}$ and $e_i^{\ast}$ are
the dual basis vectors for
$T_m^{\ast}E$ for any $m \in E$. Choose $e_1,\cdots,e_{n-1} \in
T_m(f^{-1}(t))$ so $e_1^{\ast}\wedge\cdots\wedge e_{n-1}^{\ast}$
is the volume form on $T_m(f^{-1}(t))$. Then
$e_n^{\ast}=\frac{df}{df(v(t))}$ since $v(t)$ is the unit vector
normal to $f^{-1}(t)$. Thus
\begin{eqnarray*}
Vol(E)&=&\int_E(e_1^{\ast}\wedge\cdots\wedge e_{n-1}^{\ast})\wedge e_n^{\ast}\\
&=&\int_E(e_1^{\ast}\wedge\cdots\wedge e_{n-1}^{\ast})\wedge \frac{df}{df(v(t))}
\\&=&\int_{t\in\R}dt\int_{f^{-1}(t)}\frac{(dvol)}{df(v(t))}\\
&=&\int_{\R}\frac{Vol(f^{-1}(t))}{df(v(t))}dt
\end{eqnarray*}
\hfill $\Box$

Similarly, we have the following lemma,
\begin{lm}\label{3:25}
Let $E$ be a Riemannian manifold with metric $\Gamma$,
and suppose $f:E\ra G$ is a smooth
map for which $df$ has maximal rank at generic
points in $E$.  Then
\[
Vol(E)=\int_{g\in G}\frac{Vol(f^{-1}(g))}{dg(\wedge_{j=1}^{N} df(v_j))}dg
\]
where $dg$ is the volume element given by a Riemannian metric on
$G$ and $v_1,\cdots,v_N$ are unit normal vectors to $f^{-1}(g)$ with
$N=$dim $G$.
\end{lm}

Now let us return to our particular
situation.\footnote{For other
approaches to measures on moduli spaces,
see \cite{Fine},   \cite{Fo},
\cite{Liu} and \cite{BeSe,Se1,Se2,Se3}.}
 We have a nonorientable surface $\Si$
which is the connected sum of a Riemann surface $\Si_1$ with a
nonorientable surface $\Si_2$ (either one or two copies of the
projective plane $P$). Denote
$E_1=Hom(\pi_1(\Si_1),G)$ where $\Si_1$ is a Riemann
surface of genus $\ell$ 
with one boundary component, and $E_2=Hom(\pi_1(\Si_2),G)$
where in this section $\Si_2=P$ with one disc removed. Here
$\pi_1(\Si_2)=\{x,y|x=y^2\}\cong Z$ so $Hom(\pi_1(\Si_2),G)\cong G$.
Define maps $f_i:E_i\ra G$ for $i=1,2$ by sending a
representation to its value on the loop around the boundary.
Let
$$E=\{(m_1,m_2)\in E_1\times E_2\mid f_1(m_1)=f_2(m_2)\}$$
Now we use Lemma \ref{3:25} with $E_2\cong G$
and let $f:E\ra G=E_2$ be the map $(m_1,m_2)\mapsto m_2$.
Notice that in this case the hypothesis of $f$ having maximal 
rank is valid generically 
because this reduces to $df_1: T E_1 \to \g$ being 
surjective, and it was proved by 
Goldman \cite{G2} (Proposition 3.7) that the image of 
$df_1$ at a point $A = (a_1, \dots, a_\ell, b_1, \dots, b_\ell)$ 
$\in E_1 = G^{2\ell}$ is 
$z(A)^\perp$ (the orthocomplement of the Lie algebra of the 
stabilizer of 
$A$ under the adjoint action),
and $z(A) = 0 $ for generic $A$.

We have
$$Vol(E)=\int_{G}\frac{Vol(f^{-1}(g))}{dg(\wedge_{j}df(v_j))}dg.$$
Note that the
metric on $E$ is $\pi_1^{\ast}h_1+\pi_2^{\ast}h_2$ where $h_1$ is
the metric on $E_1$ and $h_2$
the metric on $E_2$ which are given by
the equation (\ref{rpmetric}) i.e. (Haar measure)
$\times (k h(\phi))^{n/2}$. Thus we can choose  an orthonormal basis
$\{ v_j \}$  in $\g$ and $dg(\wedge_{j}df(v_j))=1$.
Thus we have $Vol(f^{-1}(g))=Vol(f_1^{-1}(g^2))Vol(f_2^{-1}(g^2))$ because
$f^{-1}(g)=f^{-1}_1(g^2)\times f^{-1}_2(g^2)$.
It follows that
\[Vol(E)=\int_{g\in G} Vol(f_1^{-1}(g^2))Vol(f_2^{-1}(g^2)) dg\]
where our moduli space is $ \M = E/G$.
Thus we have $Vol(E)=Vol(\M)Vol(G)$.

Notice that the moduli space of 
gauge equivalence classes of flat connections on a Riemann surface
with one boundary component about which the holonomy is constrained
to take a fixed value (in other words the
moduli space of parabolic bundles) is a K\"ahler manifold
(see for instance \cite{AB}) so its
symplectic volume (as specified by Witten's formula) is equal to its
Riemannian volume (for a metric derived from any metric on the Riemann
surface).
Hence our procedure will give the Riemannian volume on the
moduli space of
gauge equivalence classes of flat connections on the connected sum.

We know from \cite{w} the volume formula \cite{w}(4.114)
for a moduli space of flat
connections on a
compact orientable surface $\Si_{\ell}$ of genus $\ell$ with one
boundary component:\footnote{In Section 3.1 we give a new proof of this
formula.}
For $s \in G$, we define
\newcommand{\Rep}{{\mathcal R}}
\begin{equation} 
\M(\Si_\ell, s) = \{ (a_1, \dots, a_{2\ell})
\in G^{2 \ell} ~|~ \prod_{j = 1}^\ell a_{2j-1} a_{2j}
a_{2j-1}^{-1} a_{2j}^{-1} \in C(s) \}/G 
\end{equation}
$$ = \{ \rho \in {\rm Hom} (\pi_1(\Si_\ell - D), G) ~|~
\rho ([\partial  (\Si_\ell - D) ] ) \in C(s) \}/G $$
where $C(s)$ is the conjugacy class of $s$.
We also define
\begin{equation}
\Rep(\Si_\ell, s)  = 
\{ (a_1, \dots, a_{2\ell})
\in G^{2 \ell} ~|~ \prod_{j = 1}^\ell a_{2j-1} a_{2j}
a_{2j-1}^{-1} a_{2j}^{-1} =s \} 
\end{equation}
$$ = \{ \rho \in {\rm Hom} (\pi_1(\Si_\ell - D), G) ~|~
\rho ([\partial  (\Si_\ell - D) ] ) =  s \} $$
so that
$$\M(\Si_\ell,s) = G \Rep(\Si_\ell, s)/G. $$
It follows that \begin{equation} 
\label{e:volcomp} Vol \M(\Si_\ell, s) = Vol  \Rep(\Si_\ell, s) \frac{Vol C(s)}
{Vol G}  \end{equation}

If $ s \in G$, then Witten's formula reads
\begin{equation} \label{e:wone}
Vol(\M(\Sigma_{\ell},s))=C_1\sum_{\alpha}
\frac{1}{(\mbox{dim}\alpha)^{2\ell-1}}\chi_{\alpha}(s)\sqrt{F(s)},\end{equation}
where $\alpha$ runs over all isomorphism classes of irreducible
representations of $G$, the representation $\alpha$ has character
$\chi_{\alpha}$,  and the constant $C_1$ is
\begin{equation}\label{e:4114c}
C_1 = \frac{\sharp
Z(G)Vol(G)^{2\ell-2+1}}{(2\pi)^{\mbox{dim}
\M(\Si_\ell,s)}Vol(T)}
\end{equation} 
and \[\mbox{dim}
\M(\Si_\ell,s)=(2\ell-2)\mid G \mid + \mid G \mid - \mid T \mid,\]
where $|G|$ denotes the dimension of $G$ and
$|T|$ denotes the dimension of the maximal torus $T$.
 Note also that
 $F(s)$ is the Riemannian volume of the conjugacy class $C(s)$
through $s$ as
defined in \cite{w}(4.53)\[v=\frac{Vol(G)}{Vol(T)}v_0F(s).\]
(See Chapter 7, \cite{BGV}.) 
The Liouville volume of 
$C(s)$  is 
\begin{equation} \label{e:volconj}\sqrt{F(s)}{\rm vol} G/{\rm vol} T \end{equation}
(see \cite{AMW}, Proposition 3.6.)
Here
 $v$ represents the measure on $T/W$ obtained by pushing
down the Haar measure on $G$ (under the natural map from a group element
to its conjugacy class) and $v_0$ represents the measure on $T/W$
obtained
by restricting the metric on $\g$ to Lie($T$) and then identifying Lie($T$)
with the tangent space to $T/W$. A more detailed explanation of
$F(s)$ is given in the next section.

So (\ref{e:volcomp}) is equivalent to 
\begin{equation} \label{e:volcomp'} 
Vol \M(\Si_\ell, s) = Vol  \Rep(\Si_\ell, s) \frac{\sqrt{F(s)} }{Vol T}, 
\end{equation}
and (\ref{e:wone}) is equivalent to 
\begin{equation} \label{e:wone'}
Vol  \Rep(\Si_\ell, s)  = C_1 Vol T 
\sum_{\alpha}
\frac{1}{(\mbox{dim}\alpha)^{2\ell-1}}\chi_{\alpha}(s).
\end{equation}

Thus the volume of the moduli space of flat connections on the
connected sum of a Riemann surface with a projective plane is

\begin{equation}
Vol(\Rep(\Si))=Vol(\Rep(\Sigma_{\ell}\sharp P) )
\end{equation}
$$ = \int_{s \in G}
Vol(\Rep(\Sigma_{\ell},s))
Vol(\Rep(P,s))  \ ds 
$$

Denote by $R$ the map $G \ra G$,
\begin{equation} \label{e:rdef}
R:  g \mapsto g^2.
\end{equation}
 Weyl's integral
formula
\[
|W|\int_{G}f(g)dg=\int_{T}[F(t)\int_{G}f(gtg^{-1})dg]dt\]
gives us that \[\int_{G} f(g) dg = \int_{T/W}
F(s)f(s)\frac{Vol(G)}{Vol(T)} ds.
\]when $f$ is conjugation invariant.

Thus
\begin{eqnarray*}
Vol(\Rep(\Sigma))
&=& \int_{G}C_1 Vol T\sum_{\alpha} \frac{1}{(\dimmalpha)^{2\ell-1}}
\chi_{\alpha}(g)R_{\ast}(\sqrt{\mbox{det}\Gamma}dg)\\
&=& C_1{Vol(T)}\sum_{\alpha}
\frac{1}{(\dimmalpha)^{2\ell-1}} \int_{G} \chi_{\alpha}(g^2)
\sqrt{\mbox{det}\Gamma}dg
\end{eqnarray*} 
The above is true because $S^2-D-D$ is the double cover of $P \setminus D$; the
space of flat connections on $S^2-D-D$ (modulo based gauge transformations)
is isomorphic to $G$, so we pull back the integral to $G$ under 
the covering map.
We get \[ Vol(\Rep(\Si))
= C_1 Vol T 
 \sum_{\alpha}\frac{1}{(\dimmalpha)^{2\ell-1}} f_{\alpha}
Vol(G)(kh(\phi))^{\mid G\mid/2}
\] by using
\begin{equation}\label{e:270}
\int_{G}\chi_{\alpha}(g^2)dg=f_{\alpha} Vol(G)
\end{equation}
(cf:\cite{w} (2.70))
where $f_{\alpha} =1,-1,0$ depending on whether the representation
$\alpha$ admits a symmetric invariant bilinear form, an antisymmetric
bilinear form, or no invariant bilinear form at all.
Here $\dimmalpha$ is the dimension of the representation
$\alpha$, and its character is denoted by $\chi_\alpha$.

We get
\[Vol(\Rep(\Sigma) )
=\frac{\sharp Z(G)Vol(G)^{2\ell} }{(2\pi)^{(2\ell-1)\mid G \mid -\mid T
\mid}}H(\phi)
\sum_{\alpha}\frac{f_{\alpha}}{(\dimmalpha)^{2\ell-1}}\] where
$H(\phi)=( k h(\phi))^{\mid G\mid/2}$
and $k = \pi/\phi$. Here $h$ is a real-valued function of  the angle
$\phi$ which depends on the choice of metric on the M\"obius strip.

We compare our result
 to Witten's formula \cite{w}(4.93) \[Vol(\Rep(\Si)) =\frac{\sharp
Z(G)Vol(G)^{2\ell}}{(2\pi)^{(2\ell-1)\mid G \mid }}
\sum_{\alpha}\frac{f_{\alpha}}{(\dimmalpha)^{2\ell-1}}\]

Our formula differs from Witten's by a multiplicative factor of
$H(\phi)$ $(2\pi)^{ \mid T \mid}$. The factor
$(2\pi)^{ \mid T \mid}$ is due to Witten's choice of
a different normalization. The factor $H(\phi)$ results from our
choice of a metric on the nonorientable part of the surface (the
projective plane). Note Witten did not choose a metric on $P$.
On the other hand
our formula for $Vol(\M)$ is a function of the angle $\phi$ for a fixed
surface.

\subsection{$\Si$ is a connected sum of a Riemann surface with
two projective planes} \label{ssk}

We know from \cite{w} the volume formula for a moduli space of flat
connections on a compact orientable surface of genus $\ell$ with
two boundary components \cite{w} (4.114) (where $ s_1,s_2 \in G$
 and the  
 holonomies around the two
boundary components are fixed at $s_1$ and $s_2$): in this case
Witten's formula \cite{w} (4.114) reads as follows.
\begin{equation} \label{e:wtwo} 
Vol(\M(\Sigma_{\ell},s_1,s_2))=C_2 \sum_{\alpha}
\frac{1}{(\dimmalpha)^{2\ell}}\chi_{\alpha}(s_1)\sqrt{F(s_1)}\chi_{\alpha}(s_2)\sqrt{F(s_2)}, \end{equation}
where the constant $C_2$ is \[\frac{\sharp
Z(G)Vol(G)^{2\ell}}{(2\pi)^{\mbox{dim}
\M(\Si_\ell,s_1,s_2)}Vol(T)^2},\] and \[\mbox{dim}
\M(\Si_\ell,s_1,s_2)=(2\ell-2)\mid G \mid + 2 \mid G \mid - 2\mid T \mid.\]

Again  we define 
$$ 
\Rep(\Sigma_\ell, s_1, s_2) =  \{ \rho \in {\rm Hom}
( \pi_1 (\Sigma_\ell - D_1 - 
D_2), G) ~|~ \rho ([\partial  D_1 ] =  
  s_1 , ~ \rho ([\partial  D_2] = s_2  \}. $$
Here we have chosen  representatives 
$[\partial  D_j] $ 
for  the elements of the fundamental group represented by the loops
around the $j$-th boundary component, by connecting the basepoint 
to some point on the boundary. 
Thus
$$ \M(\Si_\ell,s_1,s_2)= (G\times G)\Rep(\Sigma_\ell, s_1, s_2)/{G\times G} $$
where $(g_1, g_2) \in G \times G$ acts on 
${\rm Hom}(\pi_1 (\Sigma_\ell- D_1- D_2), G)$ so that
if  the value of the representation of 
the loop around the $j$-th boundary component  is
$s_j$, then this value becomes
$g_j s_j g_j^{-1}$.
Following (\ref{e:volcomp}) and (\ref{e:volconj}) above, we have
\begin{equation}
{ Vol}  \M(\Si_\ell, s_1, s_2) = Vol \Rep(\Sigma_\ell, s_1, s_2) 
\frac{Vol C(s_1) Vol C(s_2)}{(Vol G)^2}
\end{equation}
$$
 = Vol \Rep(\Sigma_\ell, s_1, s_2) 
\frac{\sqrt{F(s_1)} \sqrt{F(s_2)} }{ (Vol T)^2 }. 
$$

Equivalently, Witten's formula is  
$$ Vol (\Rep(\Sigma_\ell, s_1, s_2)) = 
C_2 (Vol(T))^2 \sum_{\alpha}
\frac{1}{(\dimmalpha)^{2\ell}}\chi_{\alpha}(s_1)\chi_{\alpha}(s_2) $$

Thus the volume of the moduli space of flat connections on the
connected sum of a Riemann surface with two projective planes (equivalently,
 the connected sum with one Klein bottle) is
\begin{eqnarray*}
\lefteqn{ Vol(\Rep(\Si))=Vol(\Rep(\Sigma_{\ell}\sharp P \sharp P)) }\\
&=&\int_{G \times G} ds_1 ds_2 Vol(\Rep (\Sigma_{\ell},s_1,s_2))
Vol(\Rep (P,s_1)) Vol(\Rep (P,s_2))\end{eqnarray*}
Let $\sqrt{\mbox{det}\Gamma_1}$ and $\sqrt{\mbox{det}\Gamma_2}$ denote
the volume elements on the two copies of $P \setminus D$ respectively.
We have
\begin{equation}
Vol(\Rep(\Si))=\int_{G \times G}
Vol(\Rep(\Sigma_{\ell},g_1,g_2) )
 R_{\ast}(\sqrt{\det \Gamma_1}  dg_1)
R_\ast (\sqrt{\det \Gamma_2} dg_2)  
\end{equation}
$$
= C_2 (Vol T)^2\sum_{\alpha} \frac{1}{(\dimmalpha)^{2\ell}} \int_{G \times G}
\chi_{\alpha}(g_1)\chi_{\alpha}(g_2)R_{\ast}(
\sqrt{\mbox{det}\Gamma_1}dg_1) R_\ast (\sqrt{\mbox{det}\Gamma_2}dg_2)$$
$$= C_2 (Vol T)^2 \sum_{\alpha}
\frac{1}{(\dimmalpha)^{2\ell}} \sqrt{\mbox{det}\Gamma_1}\sqrt{\mbox{det}\Gamma_2}
(\int_{G} \chi_{\alpha}(g^2)
dg)^2 $$
where $R$ was defined by (\ref{e:rdef}).
We use equation (\ref{e:270}) and thus
the volume becomes
$$
Vol(\Rep(\Si))= \frac{\sharp Z(G)Vol(G)^{2\ell-1} }
{(2\pi)^{(2\ell)\mid G \mid -2\mid T
\mid}} \sum_{\alpha}\frac{1}{(\dimmalpha)^{2l}} f^2_{\alpha}
H(\phi_1)H(\phi_2)Vol(G)^2
$$ where
$H(\phi_i)=(k_i h(\phi_i))^{\mid
G\mid/2}$, the angles $\phi_1$ and  $\phi_2$ are the top angles of the two
triangles respectively,
 $k_i = \pi/\phi_i$,  and $f_{\alpha}$ is defined as in
Section \ref{ss:csp}.
Note that $f^2_{\alpha}=1$ if $\alpha=\bar{\alpha}$ and
$f^2_{\alpha}=0$ otherwise. We get the final formula

\[Vol(\Rep(\Si))
=\frac{\sharp Z(G)Vol(G)^{2\ell+1}}{(2\pi)^{(2\ell)\mid G \mid -2\mid T
\mid}}H(\phi_1)H(\phi_2)
\sum_{\alpha=\bar{\alpha}}\frac{1}{(\dimmalpha)^{2\ell}}.\]

We compare this with
 Witten's formula \cite{w}(4.77) \[Vol(\Rep(\Si)) =\frac{\sharp
Z(G)Vol(G)^{2\ell+1}}{(2\pi)^{(2\ell)\mid G \mid }}
\sum_{\alpha=\bar{\alpha}}\frac{1}{(\dimmalpha)^{2\ell}}\]

Our formula differs from Witten's formula by a multiplicative
factor of $H(\phi_1)$ $H(\phi_2)$
$(2\pi)^{2( \mid T \mid)}$. The factor
$(2\pi)^{2(  \mid T \mid)}$
 is due to Witten's choice
of a different normalization. The factors $H(\phi_i)$ result from
our choice of a metric on
the two projective planes.

Notice that there are two factors $H(\phi_i)$ (in contrast
 to the case when the manifold is a connected sum of a Riemann
 surface with one copy of $P$, when there is only one such
 factor).
Thus $Vol(\M)$ is a function of the angles $\phi_1,\phi_2$ for a fixed
surface $\Sigma_l$.

We may summarize our
conclusions as follows.
When we take the connected sum of a Riemann surface
with a real projective plane formed by
identifying two sides of length $L$ of a geodesic
 triangle with polar angle
$\phi$, the  volume
we have
obtained is a function of the angle $\phi$
 since the length $L$
of the side of the triangle is determined once the length of the
third edge of the
triangle is given. A pair  $(L,\phi)$ gives a conformal structure
on the triangle.
Once we fix a conformal structure, we obtain a
volume for the moduli space determined by this structure.
The situation is similar for the connected sum of a Riemann surface
with two projective planes.

\section{Witten's volume formula and Haar measure}

In this section, we will consider Witten's volume formula
for the moduli space of flat connections on a 2-manifold,
and give a new proof of these formulas using Haar measure.
 Notice that Witten did not give a 
proof exclusively using Haar measure to compute
 the volume of the moduli space.
Witten gave three arguments for the formula in \cite{w}:
\begin{itemize}
\item In Section 2, Witten used results of 
Migdal \cite{Mi} involving lattice gauge theories to 
compute the Yang-Mills partition function in two dimensions.
Migdal's result only makes sense when the Yang-Mills partition 
function has been regularized, giving a result which reduces
to the volume of the moduli space when the regularization parameter
tends to $0$.
\item In Section 3, Witten deduced the result from the 
Verlinde formula, which is a formula for the Riemann-Roch number 
of the prequantum line bundle of the 
moduli space $M$ of flat connections on a 2-manifold,
$RR(M, L^k)$.
Since the cohomology class of the symplectic form $\omega$ is the 
first Chern class of $L$, we may expect that the leading order
term in the Riemann-Roch number (as a polynomial in $k$) 
is $k^{\dim M/2} Vol(M)$. 
\item In Section 4, Witten characterized the volume of the 
moduli space using Reidemeister-Ray-Singer torsion. He used the 
fact that Reidemeister torsion is multiplicative:
If $M$ is the union of two submanifolds $M_1$ and $M_2$ with 
boundary $N$, glued along the boundary, then 
$$\tau(M) = \tau(M_1) \tau(M_2)/\tau (N)$$
where the torsion is viewed as a ratio of determinants of 
elliptic operators and the 
above quotient  makes sense in terms of the Mayer-Vietoris sequence
which computes the cohomology groups of $M$ from those of 
$M_1, $ $M_2$ and $N$.
\end{itemize}

 Since the moduli space can be identified with the space of
representations of a surface group (in other words
$\hol/G$) and a compact Lie group has a natural
Riemannian measure, the Haar measure, we can try to understand
Witten's formula for the moduli space by pushing forward the Haar
measure on $G$.\footnote{For a related treatment of the role of
Haar measure in determining volumes,
see \cite{AMM} and \cite{AMW}.}

We will use the  following  two facts  to show our argument.

$\bullet$ Weyl's integral formula (cf. \cite{bd}(4.1.11))
\[
|W|\int_{G}f(g)dg=\int_{T}[F(t)\int_{G}f(gtg^{-1})dg]dt\]
If $f$ is conjugation invariant, the formula becomes
\begin{eqnarray*}
|W|\int_Gf(g)dg&=&\int_TF(t)f(t)\frac{Vol(G)}{Vol(T)}dt\\
&=&
|W|\int_{T/W}i_{\ast}(F(t))f(i^{-1}(t))\frac{Vol(G)}{Vol(T)}i_{\ast}dt
\\&=&
|W|\int_{s\in T/W}F(i^{-1}(s))f(i^{-1}(s))\frac{Vol(G)}{Vol(T)}ds
\end{eqnarray*}
Here $i: T \ra T/W$ is the local isomorphism with
$i_{\ast}dt=:ds$.

Thus we have\[
\int_Gf(g)dg=\int_{s\in T/W}F(s)f(s)\frac{Vol(G)}{Vol(T)}ds\]

Recall that
$F(s)=F(i^{-1}(s))$ is the volume of the conjugacy class
containing $s$ as defined for example in
\cite{w}(4.53):
\[v=\frac{Vol(G)}{Vol(T)}v_0F(s),\] where $v$
represents the measure on $T/W$ obtained by pushing down the Haar measure
on $G$ (under the natural map from a group element to its
conjugacy class) and $v_0$ represents the measure on $T/W$ obtained by
restricting the metric on $\g$ to Lie($T$) (which determines a
Haar measure on $T$) and then identifying Lie($T$) with the tangent
space to $T/W$.
This is as in the picture shown below:

\begin{picture}(6,3)
\put(2.8,2){$T$} \put(5,2){$G$ }
\put(3.8,0){$T/W$}
\put(4.7,1){$u$}
\put(3,1){$i$}

\thicklines \put(3,1.7){\vector(1,-1){1}}
\put(5,1.7){\vector(-1,-1){1}}

\end{picture}

\noindent with $v=u_{\ast}dg$, and $v_0=i_{\ast}dt$,
where $dg$ and $dt$ are the Haar measures on $G$ and $T$ respectively.
Here $i:T \ra T/W$ is the local isomorphism and $u:G \ra T/W$
is the map which maps $g$ to its conjugacy class.

Explicitly we have 
\begin{equation} \label{e:Fexpl}
F(s) = {\rm det}_{\R} (1 - {\rm Ad}(s) ) 
\end{equation}
or for $s =\exp \lambda $ for $\lambda \in {\rm Lie}(T)$ 
such that $\alpha(\lambda)\ne 0 $ for any root $\alpha$, we have
$$F(\exp \lambda)  = \prod_{\alpha > 0 } 4 \sin^2 \alpha(\lambda). $$
where we assume $s \in T$ and view ${\rm Ad}(s)$ as an endomorphism of 
the orthocomplement of the Lie algebra of $T$.

\newcommand{\push}{p}

$\bullet$ Pushforward of volume under a surjective map:
Suppose
$ \push : M \rightarrow N $ is a surjective map,
  $f: N \rightarrow
R $ is a map,
and $V_{M} ~\mbox{and} ~V_{N} $ are volume forms for $M$ and $N$
respectively. Then a function $h$ characterizing the pushforward
is defined by
\[
\int_{M} f(p(x)) V_{M} = \int_{N} f(y)h(y) V_{N}, ~~\mbox{where}
~\push_{\ast}(V_{M})(y)=h(y)V_{N}. \] Roughly speaking, \(h(y)=Vol(p^{-1}(y)). \)

\subsection{$\Si$ is a one punctured Riemann surface of genus $\ell+1$}
\label{ss4.1}

The following picture shows the relations between the pushforwards.

\begin{picture}(6,6)
\put(2,5){$[G \times G]^\ell$} \put(4.5,5){~~~~$G \times G$}
\put(4,3.2){$G$~~} \put(4,1){\makebox(0,0){$T/W$}}
\put(4.7,4){$R$} \put(4.2,2){$\pi$} \put(3,4){$Q$}

\thicklines \put(3,4.7){\vector(1,-1){1}}
\put(5,4.7){\vector(-1,-1){1}} \put(4,3){\vector(0,-1){1.5}}

\end{picture}

$$
\mbox{where}~ Q(A_1,\cdots,A_\ell,B_1,\cdots,B_\ell)=\prod_{i=1}^\ell A_i
B_i A_i^{-1}B_i^{-1} $$
\begin{equation} \label{e:Rdef}
\mbox{ and} ~R(g_1,g_2)=
h g_1g_2g_1^{-1}g_2^{-1}
\end{equation}where $h \in G$ is the holonomy around the boundary.

The volume formula for a moduli space of flat connections on a
compact orientable surface with one boundary component \cite{w}(4.114)
is as follows.
For an 
element  $ s \in G$ \begin{equation}\label{w:4114}
Vol(\Rep(\Sigma_{\ell},s))=
\frac{\sharp Z(G)Vol(G)^{2\ell-1}}{(2\pi)^{(2\ell-1)\mid G\mid-\mid
T\mid} }\sum_{\alpha}
\frac{1}{(\dimmalpha)^{2\ell-1}}\chi_{\alpha}(s)\end{equation}

We want to show that Witten's formula \cite{w} (4.114) is the pushforward
of the Haar measure on $G$ i.e. prove equation (\ref{w:4114}) by
induction  on $\ell$
 (assuming that equation (\ref{w:4114}) is true for genus $\ell$, we prove it
for genus $\ell+1$).
The induction hypothesis is valid for $\ell = 0 $
because by
\cite{w} (2.43) we have
\begin{equation} \label{indstart}
\sum_\alpha (\dim \alpha) \chi_\alpha (s) = \delta (s - 1)
\end{equation}
and the  space
$\Rep(\Sigma_0,s ) $ is one point if $ s = 1$ and
empty otherwise, so its volume is
$\delta (s - 1)$.

Consider a torus $\Sigma_1$ with one boundary component.
If $\Sigma_1 - D$ is glued to $\Sigma_\ell$ along the boundary of $D$,
we obtain $\Sigma_{\ell+1}$.
(See Figure 3.)
\medskip

\medskip
\input epsf.sty
\centerline{\epsfbox{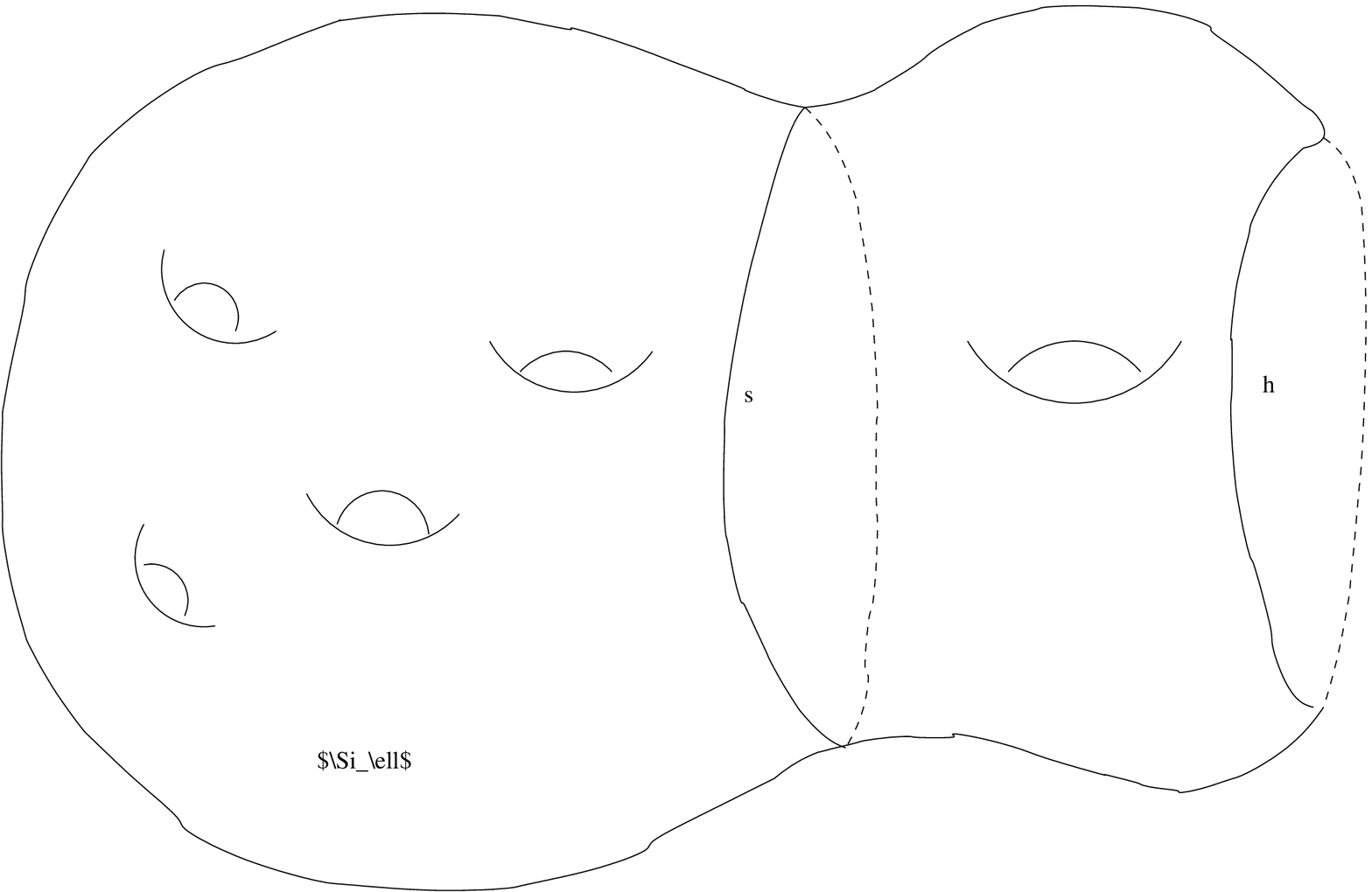}}

\noindent Figure 3: Connected sum of higher genus Riemann surface with 
 punctured torus.


Let $h \in G$ be 
the holonomy around the boundary component of
$\Si_1$ (which  becomes the boundary component of $\Sigma_{\ell + 1}$
after gluing to $\Sigma_\ell$).
We have
\begin{eqnarray*}
Vol(\Rep(\Si_{\ell+1},h))&=&Vol(\Rep(\Si_{\ell}\sharp\Si_1,h))\\
&=&\int_{s \in G}Vol(\Rep(\Si_{\ell},s))Vol(\Rep(\Si_1,s,h))ds\\
\end{eqnarray*}

Also, by definition of the pushforward we have 
\begin{equation} \label{e:pusho}
R_{\ast}(dg_1\wedge dg_2) = Vol(\Rep(\Si_1,s,h))ds.
\end{equation}

Thus we have\[ Vol(\Rep(\Si_{\ell+1},h))=\int_{G}
{Vol(\Rep(\Si_{\ell},s))}
{Vol(\Rep(\Si_1,s,h))}
 ds\]
\[ =\frac{\sharp Z(G)Vol(G)^{2\ell}}{(2\pi)^{(2\ell-1)\mid G\mid-\mid
T\mid}} \sum_{\alpha}
\frac{1}{(\dimmalpha)^{2\ell-1}} \int_{G}
\chi_{\alpha}(g)R_{\ast}(dg_1 \wedge dg_2)\]
\[=\frac{\sharp Z(G)Vol(G)^{2\ell}}{(2\pi)^{(2\ell-1)\mid G\mid-\mid
T\mid}}\sum_{\alpha} \frac{1}{(\dimmalpha)^{2\ell-1}} \int_{G\times G}
\chi_{\alpha}(hg_1g_2g_1^{-1}g_2^{-1}) dg_1 dg_2
\]
\[
=\frac{\sharp Z(G)Vol(G)^{2\ell}}{(2\pi)^{(2\ell-1)\mid G \mid
-\mid T \mid}} \sum_{\alpha}\frac{1}{(\dimmalpha)^{2\ell-1}}
\int_{G}\chi_{\alpha}(hg_1)\chi_{\alpha}(g_1^{-1}) dg_1
\frac{Vol(G)}{\dimmalpha}\] where we have
used\footnote{The equations \cite{w} (2.43), (2.48), (2.50), (2.70) and (4.50)
 are standard facts from the representation theory of
compact Lie groups, and follow from the
orthogonality relations. See also \cite{bd}.}
(cf:\cite{w}(2.50))\begin{equation}\label{e:250}
 \int_{G}\chi_{\alpha}(AuBu^{-1})du =
\frac{Vol(G)}{\dimmalpha} \chi_{\alpha}(A)\chi_{\alpha}(B)
\end{equation}
and $R$ was defined by (\ref{e:Rdef}).

Now we use (cf:\cite{w}(2.48))
\begin{equation}\label{e:248}
 \int_{G}\chi_{\alpha}(hg)\chi_{\alpha}(g^{-1})dg
 =\frac{Vol(G)}{\dimmalpha}\chi_\alpha(h)\end{equation}
to get our final formula for  $h \in G$
\begin{equation} \label{e:46a}
Vol(\Rep(\Si_{\ell+1},h))=\frac{\sharp Z(G)Vol(G)^{2\ell+1}}
{(2\pi)^{(2\ell+1)\mid G \mid -\mid T
\mid}}\sum_{\alpha}\frac{\chi_\alpha(h)}{(\dimmalpha)^{2\ell+1}} \end{equation}
We note that the final formula (\ref{e:46a}) 
is invariant under conjugation: $Vol(\M(\Si_{\ell+1},h) = 
Vol(\M(\Si_{\ell+1},khk^{-1})$ for any $k \in G$. So this formula 
only depends on the conjugacy class of $h$.

Note that Witten \cite{w} (4.114)
gives the following formula for the volume for genus $\ell+1$ with $h \in G$,
\[
Vol(\Rep(\Sigma_{\ell+1},h))= \frac{\sharp
Z(G)Vol(G)^{2\ell+1}}{(2\pi)^{(2\ell+1)\mid G\mid-\mid
T\mid}}\sum_{\alpha}
\frac{\chi_{\alpha}(h)}{(\dimmalpha)^{2\ell+1}}\]
This agrees with our result.

Our calculation shows us that Witten's formula (4.114) can be understood in terms
of
the pushforward of  Haar measure on $G$.

\subsection{$\Si$ is the connected sum of a Riemann surface of genus
$\ell$ with $P$}

We know from \cite{w} the volume formula for a moduli space of flat
connections on a compact orientable surface with one boundary
component \cite{w}(4.114), where $s \in G$ 
and the 
  holonomy around the boundary is constrained to be $s$: \[
Vol(\Rep(\Sigma_{\ell},s))=C_1 Vol T\sum_{\alpha}
\frac{1}{(\dimmalpha)^{2\ell-1}}\chi_{\alpha}(s) \]where
the constant $C_1$ is given by equation (\ref{e:4114c}).

The following picture shows the relations between the pushforwards.

\begin{picture}(6,6)
\put(2,5){$[G \times G]^\ell$} \put(4.5,5){~~~~~~~~$G$}
\put(4,3.2){$G$~~} \put(4,1){\makebox(0,0){$T/W$}}
\put(4.7,4){$R$} \put(4.2,2){$\pi$} \put(3,4){$Q$}

\thicklines \put(3,4.7){\vector(1,-1){1}}
\put(5,4.7){\vector(-1,-1){1}} \put(4,3){\vector(0,-1){1.5}}

\end{picture}

where
\begin{equation}
\label{qdef} Q(A_1,\cdots,A_\ell,B_1,\cdots,B_\ell)=\prod_{i=1}^\ell
A_i B_i A_i^{-1}B_i^{-1} \end{equation}
 and
\begin{equation} \label{rdef} R(g)=g^2. \end{equation}
By definition of the pushforward,
\begin{equation} \label{e:Fpush}
{Vol(\Rep(P,g))} dg = R_{\ast}(dg)
\end{equation}
Thus
\begin{eqnarray*}
Vol(\Rep(\Si))&=&Vol(\Rep(\Sigma_{\ell}\sharp P))
=\int_{s \in G}Vol(\Rep(\Sigma_{\ell},s)) Vol(\Rep(P,s))ds\\
&=& \int_{G} C_1 Vol T \sum_{\alpha}
\frac{1}{(\dimmalpha)^{2\ell-1}}
\chi_{\alpha}(g) R_{\ast}(dg) \\
&=& C_1 Vol T \sum_{\alpha}
\frac{1}{(\dimmalpha)^{2\ell-1}} \int_{G} \chi_{\alpha}(g^2) dg\\
&=&\frac{\sharp Z(G)Vol(G)^{2\ell}}{(2\pi)^{(2\ell-1)\mid G \mid -\mid
T \mid}} \sum_{\alpha}\frac{1}{(\dimmalpha)^{2\ell-1}} f_{\alpha}
\end{eqnarray*}
where we have used equation (\ref{e:270})
and $f_{\alpha}$ is defined in Section \ref{ss:csp}.

We get our final formula
\[Vol(\Rep(\Sigma))
=\frac{\sharp Z(G)Vol(G)^{2\ell}}{(2\pi)^{(2\ell-1)\mid G \mid -\mid T
\mid}} \sum_{\alpha}\frac{f_{\alpha}}{(\dimmalpha)^{2\ell-1}}\]

Note that Witten's formula \cite{w}(4.93) is \[Vol(\Rep(\Si)) =\frac{\sharp
Z(G)Vol(G)^{2\ell}}{(2\pi)^{(2\ell-1)\mid G \mid }}
\sum_{\alpha}\frac{f_{\alpha}}{(\dimmalpha)^{2\ell-1}}\]
Our formula differs from Witten's by a multiplicative factor of
$(2\pi)^{\mid T \mid}$ because Witten chose a
different normalization.

\subsection{$\Si$ is the connected sum of a Riemann surface of genus
$\ell$ with the Klein bottle}

Again, the following picture shows the relations between the
pushforwards. (Note that we will have a similar calculation if we use \cite{w} (4.114) with
two boundary components instead. It will be like the calculation in Section
\ref{ssk}.)

\begin{picture}(6,6)
\put(2,5){$[G \times G]^\ell$} \put(4.5,5){$G \times G$ }
\put(4,3.2){$G$} \put(4,1){\makebox(0,0){$T/W$}} \put(4.7,4){$R$}
\put(4.2,2){$\pi$} \put(3,4){$Q$}

\thicklines \put(3,4.7){\vector(1,-1){1}}
\put(5,4.7){\vector(-1,-1){1}} \put(4,3){\vector(0,-1){1.5}}

\end{picture}

\noindent where $$Q(A_1,\cdots,A_\ell,B_1,\cdots,B_\ell)=
\prod_{i=1}^\ell A_i B_i A_i^{-1} B_i^{-1}$$ and
$R(g_1,g_2)=g_1g_2g_1^{-1}g_2$.
We have that the volume of the moduli space
of flat connections on a once punctured Klein bottle with holonomy around the
puncture given by $s$ is
(by definition of the pushforward) 
\begin{equation} \label{e:pushtwo'}
{Vol(\Rep(K,g))}  dg=R_{\ast}(dg_1 \wedge dg_2).
\end{equation}
Moreover
\[ Vol(\Rep(\Si))=Vol(\Rep(\Sigma_{\ell}\sharp K))=\int_{s \in G} ds
Vol(\Rep(\Sigma_{\ell},s)) Vol(\Rep(K,s))\]
Here we have used Lemma \ref{3:25}, which is again justified 
by \cite{G2} Proposition 3.7.

Thus we have\[ Vol(\Rep(\Si))
 =\int_{G} C_1 Vol T\sum_{\alpha} \frac{1}{(\dimmalpha)^{2\ell-1}}
\chi_{\alpha}(g) R_{\ast}(dg_1 \wedge dg_2)\]
\[=C_1 Vol T \sum_{\alpha}
\frac{1}{(\dimmalpha)^{2\ell-1}} \int_{G\times G}
\chi_{\alpha}(g_1g_2g_1g_2^{-1}) dg_1 dg_2 \]
\[
=\frac{\sharp Z(G)Vol(G)^{2\ell}}{(2\pi)^{(2\ell-1)\mid G \mid
-\mid T \mid}}
\sum_{\alpha}\frac{1}{(\dimmalpha)^{2\ell-1}}
\frac{1}{\dimmalpha}
\int_{G}\chi_{\alpha}(g_1)\chi_{\alpha}(g_1) dg_1 \] where we have used
equation (\ref{e:250}).
Thus
\[Vol(\Rep(\Si)) =\frac{\sharp Z(G)Vol(G)^{2\ell+1}}
{(2\pi)^{(2\ell-1)\mid G \mid -\mid T
\mid}}\sum_{\alpha=\bar{\alpha}}\frac{1}{(\dimmalpha)^{2\ell}}\]
where we have used (cf:\cite{w}(4.50))
 \begin{equation}\label{e:450}
 \int_{G}\chi_{\alpha}(g)\chi_{\alpha}(g)dg=\left\{
\begin{array}{ll}
0 &~\mbox{if $\alpha \not=\bar{\alpha}$}\\
Vol(G)& ~\mbox{if $\alpha=\bar{\alpha}$}
\end{array}\right.
\end{equation}
Thus our final formula is
\[Vol(\Rep(\Si))=\frac{\sharp Z(G)Vol(G)^{2\ell+1} }
{(2\pi)^{(2\ell-1)\mid G \mid -\mid T
\mid}}\sum_{\alpha=\bar{\alpha}}\frac{1}{(\dimmalpha)^{2\ell}}\]

We compare our result with  Witten's formula \cite{w}(4.77): 
\[Vol(\Rep(\Si)) =\frac{\sharp
Z(G)Vol(G)^{2\ell+1}}{(2\pi)^{(2\ell)\mid G \mid }}
\sum_{\alpha=\bar{\alpha}}\frac{1}{(\dimmalpha)^{2\ell}}\] Our formula
differs from Witten's by a multiplicative factor of $(2\pi)^{\mid
G \mid+\mid T \mid} $ because Witten chose a different
normalization.

This gives us a better idea of the geometric meaning of Witten's
volume \cite{w}. It is the integral of the measure derived from the
 symplectic volume on the moduli space of flat connections on an
orientable surface and
the  pushforward volume of the Haar
measure on products of $G$.
 Since the symplectic volume of the orientable part can also
be explained as the Haar measure of the Lie group model (as we did
in section \ref{ss4.1}), this explains why the Haar measure on products
of Lie groups gives  Witten's formula for the volume on
 the moduli space of flat connections
on a nonorientable 2-manifold.

\appendix
\section{Evaluation of an integral}
In this appendix we compute the integral (\ref{e:8a}).
We assume $k$ is a positive
integer, so
 \begin{equation}\label{eq:24}
\int_{\cos(\phi/2)}^{1}\Bigl (\frac{x}{l}-\sqrt{(\frac{x}{l})^2-1}~
\Bigr )^{2k}\frac{dx}{\sqrt{1-x^2}}= S_1 + S_2, \end{equation}
where
\begin{equation} \label{eqa}
S_1 = ~\sum_{s=1}^{k} \frac{(2k)!}{(2s)!(2k-2s)!}
\int_{\cos(\phi/2)}^{1}(\frac{x}{l})^{2s}((\frac{x}{l})^2-1)^{k-s}
\frac{dx}{\sqrt{1-x^2} } \end{equation}
and
\begin{equation} \label{eqb}
S_2 = -\sum_{s=0}^{k - 1}\frac{(2k)!}{(2s+1)!(2k-2s-1)!}
\int_{\cos \phi/2}^1
(\frac{x}{l})^{2s+1}\left ( (\frac{x}{l})^{2}-1\right )^{k-s-1}\sqrt{(\frac{x}{l})^2-1}
\frac{dx}{\sqrt{1-x^2}}.
\end{equation}

We can use the trigonometric substitution $x=\cos\theta$ for
the integral (\ref{eqa}).
This yields

$$ \int_{\cos(\phi/2)}^1
 \sum_{s = 1}^k \left (\frac{x^{2s}  }{l^{2s} }\right )
 \left ( (\frac{x}{l})^2 - 1 \right )^{k-s}
\frac{dx}{\sqrt{1-x^2} } $$
$$ = \int_0^{\phi/2} \sum_{s=1}^k \frac{ (\cos^2 \theta)^s }{l^{2s} }
( \frac{\cos^2 \theta}{l^2} - 1  )^{k - s} d\theta $$
$$ =  \sum_{s = 1}^k \sum_{r = 0}^{k-s} \frac{1}{l^{2(s+r)} }
(-1)^{k-s-r} \frac{(k-s)! }{r! (k-s -r)! }
\int_0^{\phi/2}   (\cos^2 \theta)^{s+r} d\theta $$
$$ = \sum_{s = 1}^k \sum_{r = 0}^{k-s} \frac{1}{l^{2(s+r)} }
\frac{(k-s)!(-1)^{k-s-r}  }{ r! (k-s-r)!}
\Bigl [ \frac{1}{2^{s+r} }  \frac{ (2 s+2r)!}{ (s+r)! (s+r)!}  \phi/2 $$
$$
+ \frac{1}{2^{2(s+r) - 1}  }
 \sum_{p = 0}^{s+r-1} \frac{(2s+2r)!}{p! (2s+2r-2p)!}
\frac{\sin (s + r - p) \phi }{2(s+r - p)} \Bigr ] $$

The  integral (\ref{eqb}) can be obtained by the substitution
$y = \sqrt{1-x^2}$. Introducing $m^2 = 1 - l^2$ and
$y = m z$, this gives

$$
\int_{\cos(\phi/2)}^{1}\sum_{s=0}^{k-1}\frac{(2k)!}{(2s+1)!(2k-2s-1)!}
(\frac{x}{l})^{2s+1}
((\frac{x}{l})^2-1)^{k-s- 1}
 \sqrt{ (\frac{x}{l})^2-1 }
\frac{dx}{ \sqrt{1-x^2}  } $$
$$ = \sum_{s = 0 }^{k-1}\sum_{t = 0}^s  \int_{0}^{(1/m )\sin (\phi/2)}
\frac{1}{l^{2k}} \frac{s!}{t!(s-t)!} (-m^2 z^2)^t m^{2(k-s-1)} (1-z^2)^{k-s-1}
m^2 \sqrt{1-z^2}  dz $$
Now let $z = \sin \theta$, and we obtain
$$ \sum_{s = 0 }^{k-1}
\sum_{t = 0 }^s \sum_{r = 0 }^t
 \frac{m^{2(k-s+t)} }{(1-m^2)^k}  (-1)^{t+r}
 \frac{s!}{r!(t-r)! (s-t)!}
\int_0^{\sin^{-1}[(1/m )\sin (\phi/2)]} (\cos^2 \theta)^{k-s+r} d\theta $$
$$ = \sum_{s = 0 }^{k-1}
\sum_{t = 0 }^s \sum_{r = 0 }^t \frac{m^{2(k-(s-t))}}{(1-m^2)^k}  (-1)^{t+r}
 \frac{s!}{r!(t-r)! (s-t)!} \times $$
$$ \times
\Biggl [\frac{1}{2^{2(k-(s-r))} }
\ \frac{[2(k-s+r)]!}{(k-s+r)!(k-s+r)! } \theta_m
$$
$$ + \frac{1}{2^{2(k-s+r)-1} } \sum_{p = 0 }^{k- s+r -1 }
\frac{[2(k-s+r)]!}{p !(2[k-s+r] - p )! }
\frac{ \sin\{ 2 \theta_m (k - s+r- p)   \} } { 2 (k - s+r- p) }
\Biggr ]  $$
where $\theta_m = \arcsin( 1/m \sin (\phi/2) ) . $


\begin{thebibliography}{AMW}



\bibitem[AMM]{AMM} A. Alekseev, A. Malkin, and E. Meinrenken,
 {Lie group valued moment maps},
 {\em J. Differential Geom.} \textbf{48} (1998),  445--495.

\bibitem[AMW]{AMW} A. Alekseev, E. Meinrenken, and C. Woodward,
 {Duistermaat-Heckman measures and moduli spaces of flat bundles over
 surfaces},
{\em  Geom. Funct. Anal.} \textbf{12} (2002), no. 1, 1--31.


\bibitem[AB]{AB} M.F. Atiyah and R. Bott,
{ The Yang-Mills equations over Riemann surfaces}, {\em Philos.
Trans. Roy. Soc. London Ser. A} \textbf{308} (1983), no. 1505,
523--615.

\bibitem[AB2]{AB2} M.F. Atiyah and R. Bott,
{ The moment map and equivariant cohomology}, {\em Topology}
\textbf{23} (1984), no. 1, 1--28.

\bibitem[BeSe]{BeSe} C. Becker and A. Sengupta,
Sewing Yang-Mills measures and moduli  spaces over compact surfaces.
{\em J. Funct. Anal.} {\bf 52} (1998),  74--99.

\bibitem[BGV]{BGV} N. Berline, E. Getzler, M. Vergne,
{\em Heat Kernels and Dirac Operators}. Springer-Verlag,
New-York, 

\bibitem[BD]{bd} T. Br\"{o}cker and T. tom Dieck,
 {\em Representations of Compact Lie Groups},
 Graduate Texts in Mathematics, 98. Springer-Verlag, New York, 1985.

\bibitem[B]{bu} P. Buser,
{\em Geometry and Spectra of Compact Riemann Surfaces},
Birkh\"auser
(Progress in Math, Vol. 106),  1992.

\bibitem[Cox]{cox} H.S.M. Coxeter, {\em Introduction to Geometry}, Wiley, New York,
1989.

\bibitem[dC]{dC} M. do Carmo, {\em Differential Geometry of
Curves and Surfaces}, Prentice-Hall, 1976.

\bibitem[d'H]{d'H} E. d'Hoker, {\em String Theory}, in 
{\em Quantum Fields and Strings: A Course for 
Mathematicians} (ed. P. Deligne {\em et al.}), Vol. 2, AMS, 1999, p. 807-1012.


\bibitem[Fine]{Fine}
D. Fine,
Quantum Yang-Mills on the two-sphere.
{\em Comm. Math. Phys.}  {\bf 134} (1990),  273--292.




\bibitem[Fo]{Fo}
R. Forman,
Small volume limits of $2$-d Yang-Mills.
{\em Comm. Math. Phys.} {\bf  151}  (1993),  39--52.

\bibitem[Fr]{Fr} D. Freed,
Reidemeister torsion, spectral sequences,
and Brieskorn spheres. {\em J. Reine Angew. Math.} {\bf 429} (1992), 75--89.

\bibitem[Go1]{G2}
W.M. Goldman, {The symplectic nature of fundamental groups of
surfaces}, {\em Adv. Math.} \textbf{54} (1984), 200-225.



\bibitem[GM]{GM} W.M. Goldman and J.J. Millson,
{The deformation theory of representations of fundamental
groups of compact K\"{a}hler manifolds}, {\em Bull. Amer. Math. Soc.}
(N.S.) \textbf{18} (1988), no. 2, 153--158.








\bibitem[Liu]{Liu} K. Liu,
 Heat kernel and moduli space,  {\em Math. Res. Lett.}
{\bf  3} (1996),  743--762.


\bibitem[Ma]{ma} W. S. Massey,
 {\em Algebraic Topology: An Introduction},
 Graduate Texts in Mathematics, 56. Springer-Verlag, New York, 1967.

\bibitem[Mi]{Mi} A.A. Migdal, Recursion equations in gauge field
theories. {\em Sov. Phys. JETP} {\bf 42}, 413-418 (1976).

\bibitem[Se1]{Se1} A. Sengupta, Yang-Mills on surfaces
with boundary: quantum theory and symplectic limit.
{\em Commun. Math. Phys.} {\bf 183} (1997) 661-705.

\bibitem[Se2]{Se2} A. Sengupta, Sewing symplectic volumes for flat connections
over compact surfaces. {\em J. Geom. Phys.} {\bf 32} (2000) 269-292.

\bibitem[Se3]{Se3} A. Sengupta, The Yang-Mills measure and
symplectic structure over spaces of connections.
{\em Quantization of singular symplectic quotients}, 329-355,
Birkh\"auser (Progress in Mathematics vol. 198), 2001.



\bibitem[Weis]{Weis} Eric W. Weisstein, Spherical trigonometry. From 
{\em MathWorld}-- A Wolfram Web Resource.
\begin{verbatim}
http://mathworld.wolfram.com/SphericalTrigonometry.html
\end{verbatim}

\bibitem[W]{w} E. Witten, {On quantum gauge theories in
two dimensions}, {\em Commun. Math. Phys.}  \textbf{141} (1991), 153--209

\end{thebibliography}
\end{document}